\documentclass[12pt, reqno]{amsart}
\usepackage{amsmath, amsthm, amscd, amsfonts, amssymb, graphicx, xcolor}
\usepackage{tikz}
\usepackage[bookmarksnumbered, colorlinks, plainpages]{hyperref}
\usepackage{tikz}
\usepackage{enumitem}
\usepackage{bm}
\usepackage{float}
\newtheorem{theorem}{Theorem}[section]
\newtheorem{lemma}[theorem]{Lemma}
\newtheorem{proposition}[theorem]{Proposition}
\newtheorem{corollary}[theorem]{Corollary}
\theoremstyle{definition}
\newtheorem{definition}[theorem]{Definition}

\theoremstyle{remark}

\numberwithin{equation}{section}
\newtheorem{notation}{Notation}[section]

\begin{document}

\setcounter{page}{1}

%

\centerline{}

\centerline{}


\title[]{Maximum number of  cells in a matroid subdivision of hypersimplex }

\author[]{ Mahyar Hasani}

\footnote{This work was partially supported by INSF under Grant No. 4049675.}


\address{$^{3}$ Department of Mathematics, IPM - Institute for Research in Fundamental Sciences, Tehran, Iran}
\email{\textcolor[rgb]{0.00,0.00,0.84}{mahyarhasani@ipm.ir}}

%
%

\begin{abstract}

David Speyer proved that
the maximum number of  cells in a  regular matroid subdivision of the hypersimplex
$\Delta(m,t)$
is
$\binom{t-2}{m-1}$
, using the language of tropical geometry and properties of the Tutte polynomial of a matroid \cite[Theorem 6.1]{Speyer}.
In this paper, we prove this result in general, i.e. for arbitrary matroid subdivisions including the non-regular ones. We perform  a different  combinatorial way using our results in \cite{HHH} and  independent of tropical language. 
\newline
\end{abstract} \maketitle


\tableofcontents

\section{Introduction}

Triangulations and subdivisions are among  the most important notions in convex geometry as well as  in algebraic geometry and topology. For example in toric geometry,  we can associate an algebraic  variety to a triangulation of a special polytope.   More generally we can consider face-fitting decompositions, in which we do not need the 
smaller polytopes to be simplices anymore. One of the most important special cases is a subdivision of hypersimplex 

\begin{center}
	$\Delta(m,t)= 
	\left\{(x_{1},x_{2},...,x_{t})\in \mathbb{R}^t:x_{1}+...+x_{t}=m, 0\le x_{i}\le 1\right\}$  
\end{center}

by matroid polytopes. In fact, the hypersimplex is itself the polytope of the uniform matroid $U_{m,t}$ and its subdivision is also considered in matroid theory\cite{oriented}.
On the other hand one can construct toric varieties from 
matroid polytopes
involved in a matroid subdivision and  attach them to get the so called stable varieties or "shas"(stable hyperplane arrangements) as explained in \cite{Alex}. A natural question is what is the maximum number of irreducible components in a "shas" obtained by a matroid subdivision of hypersimplex 
$\Delta(m,t)$?
Or, in another words, what 
is the maximum number of polytopes in a matroid subdivision of the hypersimplex? We will recall some definitions from matroid theory:

\begin{definition}
	A \emph{matroid} is a pair 
	$P=(\left\{1,...,t\right\},\varGamma)$,
	where 
	$\varGamma$
	is  a collection  
	of subsets of
	$\left\{1,...,t\right\}$
	satisfying the properties below:

	\begin{enumerate}
		\item 
		$ \emptyset \in \varGamma $
		\item 
		If $I\in \varGamma$ and $J\subseteq I$, then $J\in \varGamma$
		\item 
		If $I,J \in \varGamma$  and $|I| < |J|$ , there exists $j\in J-I$
		such that $I\cup  \left\{j\right\} \in  \varGamma $.

	\end{enumerate}
	
	We also denote the matroid simply by $\Gamma$.
	
\end{definition}

\begin{definition}
	If the  pair 
	$P=(\left\{1,...,t\right\},\varGamma)$ is a matroid, 
	the elements of $\varGamma$ are  called independent sets of the matroid.
	
\end{definition}

Matroids can also be defined via their bases:

\begin{definition}
	Let the  pair 
	$P=(\left\{1,...,t\right\},\varGamma)$ be a matroid, a subset 
	$I\subseteq \left\{1,...,t\right\}$ is called a \emph{basis} if it is a  
	maximal independent set. We denote the set of all bases by $\beta(P)$ or $\beta(\Gamma)$.
	
\end{definition}




\begin{definition}
	The \emph{matroid polytope} of a matroid $P=(\left\{1,...,t\right\},\varGamma)$ is the convex hull of the points $e_I$ with $I\in \beta(P)$. It is a polytope in $\mathbb{R}^t$ and is denoted by $M_P$ or $M_\Gamma$. 
\end{definition}

Two operations can be performed on matroids, one is \emph{deletion} and the other is \emph{contraction} :

\begin{definition}\label{del}
	Let $P=(\left\{1,...,t\right\},\varGamma)$ be a matroid and $X\subseteq \left\{1,...,t\right\}$. Put 
	\begin{center}
		$\varGamma|X=\left\{ I\subseteq X: I\in \varGamma\right\}$
	\end{center}
	
	Then the pair 
	$(X,\varGamma|X)$ is also a matroid, called the \emph{deletion of}$X^c$  and is denoted by $P|X$.
	
\end{definition}

\begin{definition}\label{cont}
	Let $P=(\left\{1,...,t\right\},\varGamma)$ be a matroid and $X\subseteq \left\{1,...,t\right\}$.
	
	Put 
	
	\[
	\varGamma/X = \left\{ I \subseteq X^c: P|X \text{ has a basis } B \text{ such that } B \cup I \in \varGamma \right\}
	\]
	
	Then the pair 
	$(X^c,\varGamma/X)$ is also a matroid, called the \emph{contraction} of $X$  and is denoted by $P/X$.

\end{definition}

The notion of  rank function for a matroid is essential in our study:

\begin{definition}
	Let $P=(\left\{1,...,t\right\},\varGamma)$ be a matroid and $X\subseteq \left\{1,...,t\right\}$. The rank of $X$, denoted by $r_{P}(X)$, is defined to be the cardinality of a basis of $P|X$.
\end{definition}

We will work with two kinds of matroids:

\begin{definition}
		A  matroid $P=(\left\{1,...,t\right\},\varGamma)$ of rank $m$ is \emph{linear} if there are
		  vectors 
		 	$f_1,...,f_t$ in $\mathbb{C}^m$
		 	such that  a subset
		 	$I\subseteq \left\{1,...,t\right\}$
		 	is independent if and only if
		 	the set $\left\{f_i:i\in I\right\}$	is linearly independent over $\mathbb{C}$.
\end{definition}
 There is a notion of connectedness in matroid theory which is useful for our aim.

\begin{definition}
	A  matroid $P=(\left\{1,...,t\right\},\varGamma)$  is
	connected if one cannot decompose the ground set 
\begin{center}
		$\left\{1,...,t\right\}=E_1\cup E_2$
\end{center}
	such that
	 a subset
	$I\subseteq \left\{1,...,t\right\}$
	is independent if and only if
	
	\begin{center}
		$I=I_1\cup I_2$
	\end{center}
	where
	$I_1\subseteq E_1$
	,
		$I_2\subseteq E_2$
	 and $I_1$, $I_2$
	are independent.

\end{definition}


\begin{definition}
	A collection of maximal matroid polytopes $\left\{Q_j: j\in J\right\}$ in $\mathbb{R}^t$ is a \emph{matroid decomposition} of $\Delta(m,t)$
	if $\bigcup\limits_{j\in J} Q_j=\Delta(m,t)$ and they intersect face-fittingly; i.e., for each  pair  $Q_i$ and $Q_j$ , $Q_i\cap Q_j$ is either empty or a face of both $Q_i$ and $Q_j$.
\end{definition}

It is proven by 
Herrmann,  Jensen, Joswig, and Sturmfels
that 
the maximum number of  cells in a matroid subdivision of the hypersimplex
$\Delta(3,t)$ is   not bounded by a constant \cite{Tropical}.  Valery Alexeev also proved that 
this number is bounded by $(t-3)^2$  for $\Delta(3,t)$ \cite{Alex}. 
Finally, David Speyer proved that for all $m,t$   
the maximum number of  cells in a regular matroid subdivision of 
$\Delta(m,t)$
is
$\binom{t-2}{m-1}$.
 In this article we extend this
result to all arbitrary(including non-regular)  matroid subdivisions of 
$\Delta(m,t)$ via the ideas written in 
\cite{HHH}. We hope our new approach also helps us to answer a question about  the Dressian fan  in tropical geometry, denoted by $Dr(m,t)$. We briefly state the question below:

 Every element in the Dressian fan
 $Dr(m,t)$ determines a matroid subdivision for
 $\Delta(m,t)$. The number, \emph{spread} of an element in the Dressian, is the 
 number of maximal cells in the corresponding matroid subdivision of hypersimplex. Question 35 in \cite{Dressian} states:
 \begin{itemize}
 	\item 
 What is the maximal spread of a ray of $Dr(m,t)$ (at least for m = 3)?
   \end{itemize}

  In the question above, rays  correspond to  subdivisions that
  cannot be coarsened, i.e. 
  if you want to coarsen it, you will only get the whole hypersimplex $\Delta(m,t)$.
  We hope that the recursive process that we used for the previous problem can also be used to solve the problem above.

\section{Relation to inverse limit of point configurations}

The notion of inverse limit for  
an action of an algebraic group on an algebraic variety was studied in \cite{Th99} using Geometric Invariant Theory \cite{Mumford}
and its variation \cite{Th94,DH98}. For a reference concerning this matter  and the special case of point configurations see
\cite{BH}.

\begin{definition}
	A point configuration with parameters $m,t$ is a tuple $(p_1,...,p_t)$ such that 
	$p_i\in \mathbb{P}^{m-1}$ for every $i=1,...,t$. We also assume that we have at least $m+1$ general points in this configuration.
\end{definition}

One can consider finitely many quotients of the natural action of $PGL(m)$ on $(\mathbb{P}^{m-1})^t$, provided by 
Geometric Invariant Theory, and consider their inverse limit.  In fact, for every weight 
$\theta \in \Delta(m,t)$, we have a corresponding moduli space 
$M^{\theta}$ parametrizing $\theta$-semistable representations up to S-equivalence. Based on the Geometric Invariant Theory literature, they form an inverse system and we denote this inverse limit by 
$\overline{M}(m,t)$.

 In fact, the walls 
$W^l_I$ for 
$I\subseteq \left\{1,...,t\right\}$ , $1\leq l<m$  decompose the hypersimplex 
$\Delta(m,t)$ into finitely many chambers. For each chamber in this decomposition, we have a distinct moduli space.
We can specify a representative $\theta$ for each chamber; this representative determines the corresponding moduli space $M^\theta$.
The inverse limit
$\overline{M}(m,t)$
is a closed subvariety of
$\prod\limits_{\theta\in C(m,t)}  M^\theta$,
where 
$C(m,t)$  is the set which has exactly one weight as a representative from each chamber. In this way, we can denote each element of the inverse limit 
$\overline{M}(m,t)$
as 
$(\varphi_{\theta})_{\theta\in C(m,t)}$ such that $\varphi_{\theta}\in M^{\theta}$.

We need to establish some notations:

\begin{notation}

	Let $e_1,...,e_t$ be the standard basis vectors of $\mathbb{R}^t$. Suppose $I\subseteq \left\{1,...,t\right\}$, $0\leq l<m$ and $\theta=(\theta_1,...,\theta_t) \in \mathbb{R}^t$.
	\begin{center}
		\begin{enumerate}
			\item 
			$e_I=\sum\limits_{i\in I} e_{i}$

			\item 
			$\theta(I)=\theta_I=\sum\limits_{i\in I} \theta_{i}$ 
			\item 
			$W^l_I=\left\{(\theta_1,...,\theta_t)\in \mathbb{R}^t : \sum\limits_{i\in I} \theta_{i}=l \right\}$ 
			

		\end{enumerate}
	\end{center}
\end{notation}

In the above notations, the  particular case  $I=\left\{i \right\}$ is of interest to us. In this case  we have

\begin{center}
	\begin{itemize}
		\item 
			$W^0_{\left\{i\right\}} \cong \Delta(m,t-1)$
			\item 
		$W^1_{\left\{i\right\}} \cong \Delta(m-1,t-1)$
		\end{itemize}
\end{center}
We will work with both of them multiple times.

  In the paper 
\cite{HHH}, we construct algebraic morphisms

\begin{center}
	$
	\overline{M}(m, t) \longrightarrow \overline{M}(m, t-1)
	$

	$
	\overline{M}(m, t) \longrightarrow \overline{M}(m-1, t-1)
	$

\end{center}

for every 
$1 \leq i\leq t$
which are called
 \emph{deletion} and \emph{contraction} maps. They are
 analogous 
 to the well-known
 maps 
 of the
  Kapranov-Chow quotients of Grassmannians \cite{Kap1}. One can combine the 
   \emph{deletion} and \emph{contraction} maps
   for a fixed 
   $1 \leq i\leq t$
    simultaneously
    to get the natural morphism
   \begin{center}
   	$
   	\overline{M}(m, t) \longrightarrow \overline{M}(m, t-1)\times \overline{M}(m-1, t-1)
   	$.
   	\end{center}

 If we have a point of the
  inverse limit
  	$
  \overline{M}(m, t)$, this gives us a matroid decomposition of the hypersimplex
  $\Delta(m,t)$  for $m=2,3$ by \cite[Theorem 1.14]{HHH}, we can apply the morphisms above to 
  get  elements in 
  	$
  \overline{M}(m, t-1)$
  and 
  	$
  \overline{M}(m-1, t-1)$ and decompositions for
  $\Delta(m,t-1)$ and $\Delta(m-1,t-1)$ respectively. 
  
  In fact these morphisms
  are derived from natural maps:
  
  \begin{center}
  	
  	$\eta_i:\prod\limits_{\theta\in C(m,t)}  M^\theta\longrightarrow \prod\limits_{\theta\in C(m,t-1)}  M^\theta$
  	
  	$
  	(\varphi_{\theta})_{} \mapsto (\varphi_{\theta}^{D_i})_{}
  	$
  			
  	\end{center}
  	
  	\begin{center}
  
  	$\gamma_i:\prod\limits_{\theta\in C(m,t)}  M^\theta\longrightarrow \prod\limits_{\theta\in C(m-1,t-1)}  M^\theta$
  	
  	$
  	(\varphi_{\theta})_{} \mapsto (\varphi_{\theta}^{C_i})_{}
  	$
  \end{center}
  where
  $\varphi_{\theta}^{D_i}$ and
  $\varphi_{\theta}^{C_i}$
  are obtained from the representation 
  $\varphi_{\theta}$ by applying the 
  \emph{deletion} and \emph{contraction} process over point $i$
   for each
  $1 \leq i\leq t$  as explained in 
  \cite[Definition
  3.14,Definition 3.15]{HHH}. 
   In this definition, 
  chambers of the hypersimplices
  $\Delta(m,t-1)$ and $\Delta(m-1,t-1)$
  are obtained with intersecting of chambers of 
  $\Delta(m,t)$ with 
  $W^0_{\left\{i\right\}}$,$W^1_{\left\{i\right\}}$ respectively. 
  So if we have a point of the inverse limit   
  $(\varphi_{\theta})_{\theta \in C(m,t)}$
  in
  $
  \overline{M}(m, t)$
   the images would be tuples of the form
   $(\varphi_{\theta}^{D_i})_{\theta\in C(m,t-1)}$
   and
   $(\varphi_{\theta}^{C_i})_{\theta\in C(m-1,t-1)}$ respectively.

  For 
  a  point inside the  inverse limit
  $\overline{M}(m,t)$,
   we denote the
  maximum number of cells in the corresponding  subdivision of hypersimplex
  $\Delta(m,t)$ by K(m,t). Then we can get a recursive formula if we prove the items below:

  \begin{enumerate}
  	\item 
  	Every component of 
  	 $(\varphi_{\theta}^{D_i})_{\theta\in C(m,t-1)}$
  	and
  	$(\varphi_{\theta}^{C_i})_{\theta\in C(m-1,t-1)}$ is obtained from one component of 
  	 $(\varphi_{\theta})_{\theta \in C(m,t)}$.

  	\item 
  		Every  component of the 
  		initial element of
  		$(\varphi_{\theta})_{\theta \in C(m,t)}$
  		in $\overline{M}(m,t)$ appears in at least one(or both) of the components of the images
  		$(\varphi_{\theta}^{D_i})_{\theta\in C(m,t-1)}$
  		and
  		$(\varphi_{\theta}^{C_i})_{\theta\in C(m-1,t-1)}$
  		in
  	$\overline{M}(m,t-1)$ and 
  	$\overline{M}(m-1,t-1)$ respectively.
  \end{enumerate}

  If we prove these two facts, then we have:
  
  \begin{center}
  	
  	\begin{itemize}
  		\item 
  		$K(m,t-1)+K(m-1,t-1)\geq K(m,t)$
  		
  	\end{itemize}
  	  \end{center}
  	
  	To see why this is true, take a  tuple with 
  	$K(m,t)$ maximal cells
  	 which is  a point of the inverse limit
  	$\overline{M}(m,t)$, the items (1), (2) in the previous paragraph show that 
  	we can decompose the cells of this tuple into union of two sets $A\cup B$ such that $A$  and $B$ are the
  	maximal cells of the images in 
  	$\overline{M}(m,t-1)$
  	and
  	$\overline{M}(m-1,t-1)$
  	respectively, so
  	
  	\begin{center}
  		
  		$K(m,t)=|A\cup B|\leq |A|+|B|$
  	\end{center}

  	 Thus using induction on $m+t$ and Pascal identity, we get
  	 
  	 \begin{center}
  	 	
  	 $K(m,t)\leq \binom{t-3}{m-1}+ \binom{t-3}{m-2}=\binom{t-2}{m-1}$
  	
  	\end{center}

   for 
   $m\geq t+1$ and $m\geq 2$ and the bound is achieved. The base of induction is also clear ($m=2,t=4$).

  It is quite easy to see that both of (1), (2) are satisfied for point configurations or linear matroids. 
  With this motivation we prove (1),(2) for general matroids in theorem \ref{main} and state this bound for every matroid decomposition of hypersimplex
  $\Delta(m,t)$
  independent of inverse limit  and regularity concept.


  \section{Proof that the bound is working }
  We first  establish some notation and lemmas, specially the lemma below, which is well-known in matroid theory.
  \begin{lemma}\label{connecteddim}
  	Let
  	$(\left\{1,...,t\right\},\varGamma)$
  	be an abstract matroid,
  	the matroid polytope
  	$M_{\Gamma}$
  	is of maximal dimension if and only if the matroid $\Gamma$ is connected.
  \end{lemma}

  \begin{proof}
  	Dimension of the matroid polytope can be computed using \cite[Proposition 2.4]{DMP},i.e
  	\begin{center}
  		$\dim(M_{\Gamma})=t-C(\Gamma)$
  	\end{center}
  	where
  	$C(\Gamma)$
  	is the number of connected components
  	of the matroid
  	$\Gamma$. Notice that a matroid  $\Gamma$ is connected if and only  of $C(\Gamma)=1$. Since the matroid polytope $M_{\Gamma}$ is a sub-polytope of the hypersimplex $\Delta(m,t)$ and $\dim(\Delta(m,t))=t-1$ the result is driven.
  \end{proof}

  \begin{notation}
  		If 
  		$P=(\left\{1,...,t\right\},\varGamma)$ is a matroid and $1\leq i\leq t$, we can operate the 
  	 \emph{Deletion} and \emph{Contraction} process , over point $i$ and get two new matroids
  	  $\Gamma^{D_i}$ and $\Gamma^{C_i}$
  	  with indices
  	  $\left\{1,2,..., i-1,i+1,...t\right\}$. 
  	 
  \end{notation}
  
  %
  	%
  %
  %
  
  We also need  the definition of isomorphism between matroids:
  
  \begin{definition}\label{isomat}
  	Two abstract matroids 
  	$(\left\{1,...,t\right\},\varGamma)$ and 
  	$(\left\{1,...,t\right\},\varGamma^\prime)$
  	are isomorphic if there is a bijective map 
  	$f:\left\{1,...,t\right\}\to \left\{1,...,t\right\}$
  	such that 
  	for every subset
  	$I\subseteq\left\{1,...,t\right\}$, 
  	$I\in \varGamma$
  	if and only if 
  	$f(I)\in \varGamma^\prime$.
  \end{definition}
  
  \begin{definition}

  	We say two matroid subdivisions
  	
  	\begin{center}
  	$\left\{M_{\Gamma_1},...,M_{\Gamma_q}\right\}$	
  		
  	\end{center}

  \begin{center}
  	,
  \end{center}

  	\begin{center}
  		$\left\{M_{\Gamma^\prime_1},...,M_{\Gamma^\prime_q}\right\}$	
  		
  	\end{center}

  	 are isomorphic if there is a bijection between their cells, i.e
  	 there are isomorphisms of matroids
  	 	$f_i:\Gamma_i\to \Gamma_i^\prime$
  	  for
  	  $1\leq i\leq q$.

  \end{definition}
  
  Let $D(m,t)$ be the set of all matroid subdivisions of the hypersimplex $\Delta(m,t)$ up to the notion of isomorphism defined in the above. We will construct the map  below 
  analogous the inverse limit one:

  \begin{center}
  	\[
  D(m, t)  \longrightarrow  D(m, t-1)\times D(m-1, t-1)
  	\]
  	
  	\[
  	(M_{\Gamma_s})^q_{s=1}
  	\longmapsto
  (	(M_{\Gamma_s}\cap \Delta(m,t-1))^q_{s=1}, 	(M_{\Gamma_s}\cap \Delta(m-1,t-1))^q_{s=1})
  	\]
  \end{center}

  Here we should notice that since the decomposition 
  	$(M_{\Gamma_s})^q_{s=1}\in D(m,t)$
  is face-fitting, after intersecting with
  $\Delta(m,t-1)$ and $\Delta(m-1,t-1)$, it will remain face-fitting.
  We should also prove that both of the components 
  $M_{\Gamma_s}\cap \Delta(m,t-1)$ and $M_{\Gamma_s}\cap \Delta(m-1,t-1)$ are matroid polytopes. We will prove this in the proposition below:

  %

  \begin{proposition}\label{intersection}
  	
  	If 
  	$(\left\{1,...,t\right\},\varGamma)$
  	 is an abstract matroid  such that $M_{\Gamma}$ is of maximal dimension and $1\leq i\leq t$,
  	then
  	
  	\begin{itemize}
  		\item 
  		$M_{\Gamma}\cap W^0_{\left\{i\right\}}\cong M_{\Gamma^{D_i}}$
  		\item 
  		$M_{\Gamma}\cap W^1_{\left\{i\right\}}\cong M_{\Gamma^{C_i}}$
  		
  	\end{itemize}
  	
  \end{proposition}

\begin{proof}
	Defenition of matroid polytope 
	$M_{\Gamma}$
	says:
	\begin{center}
			$M_{\Gamma}\cap W^0_{\left\{i\right\}}=conv\left\{e_I:I\in \beta(\Gamma), i\notin I \right\}$
			
				$M_{\Gamma}\cap W^1_{\left\{i\right\}}=conv\left\{e_I:I\in \beta(\Gamma), i\in I \right\}.$
			
	\end{center}
	Based on the assumption the matroid polytope $M_{\Gamma}$ is of maximal dimension. Thus lemma \ref{connecteddim} implies that 
	the matroid
	$\Gamma$ is connected and the element $i$ is neither loop nor coloop. So we have
		\begin{center}
		$\beta(\Gamma^{D_i})=\left\{I\in \beta(\Gamma): i\notin I \right\}$
		
		$\beta(\Gamma^{C_i})=\left\{I-\left\{i\right\}:I\in \beta(\Gamma), i\in I \right\}$
		
	\end{center}
	 Now it is easy to see the isomporphisms.
\end{proof}
 Now if we prove
items (1), (2) of the theorem below , then we get the required bound by the same way as before.
 \begin{theorem}\label{main}
	With the notations above, we have:
	
	\begin{enumerate}
		\item 
		Every  maximal component of
		$(M_{\Gamma_s}\cap \Delta(m,t-1))^q_{s=1}$
		and
		$(M_{\Gamma_s}\cap \Delta(m-1,t-1))^q_{s=1}$
		 is obtained from one component of 
		$(M_{\Gamma_s})^q_{s=1}$.

		\item 
		Every   component of the 
		initial decomposition 
			$(M_{\Gamma_s})^q_{s=1}$
		in $D(m,t)$ appears in at least one of the maximal components of the images
			$(M_{\Gamma_s}\cap \Delta(m,t-1))^q_{s=1}$
		and
		$(M_{\Gamma_s}\cap \Delta(m-1,t-1))^q_{s=1}$
		in
		$D(m,t-1)$ and 
		$D(m-1,t-1)$ respectively.
	\end{enumerate}
\end{theorem}

\begin{proof}
	
	Item (1) is clear. For (2) notice that   the intersection 
	$M_{\Gamma_s}\cap \Delta(m,t-1)$
	and
	$M_{\Gamma_s}\cap \Delta(m-1,t-1))$
	could  be non-maximal, in this case they are not  computed as a maximal component. We should prove  that for every component 
	$M_{\Gamma_s}$ at least one of 
	$M_{\Gamma_s}\cap \Delta(m,t-1)$
	or
	$M_{\Gamma_s}\cap \Delta(m-1,t-1)$
	is maximal. This is proven in the next proposition.
	
	\end{proof}

  \begin{proposition}
  	At least one of the polytopes
  	 $M_{\Gamma_s}\cap \Delta(m,t-1)$
  	or
  	$M_{\Gamma_s}\cap \Delta(m-1,t-1)$
  	is maximal dimensional.
  \end{proposition}
  
  \begin{proof}
  	Based of proposition \ref{intersection} we have to prove at least one of the polytopes
  	$M_{\Gamma_s^{D_i}}$ or 
  	$M_{\Gamma_s^{C_i}}$
  	is of maximal dimensional. By lemma \ref{connecteddim} we should prove that 
  	either
  		$\Gamma_s^{D_i}$ or 
  	$\Gamma_s^{C_i}$
  	 is connected. Since
  	 	$(M_{\Gamma_s})^q_{s=1}\in D(m,t)$,
  	  the matroid polytope $M_{\Gamma_s}$ is of maximal dimension, so 
  	$\Gamma_s$ is connected. As known in matroid theory, Theorem 4.3.1 in Oxley book \cite{Ox} implies that either 
  	  $\Gamma^{D_i}$ or 
  	  $\Gamma^{C_i}$
  	  is connected.
  \end{proof}

  \begin{corollary}
  Maximum number of maximal cells in an arbitrary matroid decomposition of hypersimplex $\Delta$(m,t) is
  $\binom{t-2}{m-1}$.
  	
  \end{corollary}

%

{\bf Acknowledgement.} 

I would like to express my gratitude to my advisors, Prof. Lutz Hille and Prof. S. Mohammad Hadi Hedayatzadeh, for their constant support and invaluable guidance. In particular, I am grateful to Prof. Lutz Hille for introducing me to this interesting problem and the concept of inverse limit, which played a key role in shaping my research.

\bibliographystyle{amsplain}

\end{document}